\documentclass[12pt]{article}
\usepackage{amsfonts}

\usepackage{latexsym}
\usepackage{amssymb}
\usepackage{epsfig}
\usepackage{amsmath}
\usepackage{amsthm}
\usepackage[mathscr]{eucal}
\usepackage{subfigure}
\usepackage{graphicx}
\usepackage{hyperref}

\newtheorem{Lemma}{Lemma}

\newtheorem{Theorem}[Lemma]{Theorem}

\newtheorem{Corollary}[Lemma]{Corollary}
\newtheorem{Definition}{Definition}

\renewcommand{\qed}{\hfill{\ \ \rule{2mm}{2mm}} \vspace{0.2in}}

\begin{document}

\title{Graph extensions, edit number and regular graphs}
\author{ \textbf{Ghurumuruhan Ganesan}
\thanks{E-Mail: \texttt{gganesan82@gmail.com} } \\
\ \\
New York University, Abu Dhabi }
\date{}
\maketitle

\begin{abstract}


A graph~\(G\) on~\(n\) vertices is said to be extendable if~\(G\) can be modified to form a new graph~\(H\) on more than~\(n\) vertices, while preserving the degrees of the vertices common to~\(G\) and~\(H.\) The added vertices all have the same degree and we define edit numbers to quantify the amount of modification needed to obtain the extended graph. Characterizing graphs with least possible edit numbers, we obtain that graphs with zero edit number can be extended using regular graphs. We also describe an iterative algorithm to construct connected regular graphs on arbitrarily large vertex sets, starting from the complete graph on a fixed set of vertices.

\vspace{0.1in} \noindent \textbf{Key words:} Graph extensions, edit number of graphs, connected regular graphs.

\vspace{0.1in} \noindent \textbf{AMS 2000 Subject Classification:} Primary:
60J10, 60K35; Secondary: 60C05, 62E10, 90B15, 91D30.
\end{abstract}

\bigskip

\setcounter{equation}{0}
\renewcommand\theequation{\thesection.\arabic{equation}}
\section{Introduction}
The process of transforming one graph into another with additional
vertices and satisfying desirable properties has been studied
before in many contexts. Erd\H{o}s and Kelly (1963), Akiyama et al (1985)
estimated the minimum number of additional vertices needed to
convert a given graph into a regular graph. Bodlaender et al (2003)
determine the time needed to construct such a regular supergraph.

In many applications, it is also important to obtain the transformed graph
with as few edit operations as possible, typically comprising vertex addition/deletion and/or edge addition/deletion.
The degree sequence of the resulting graph is often very different from the original graph.
In this paper, we consider a slightly different problem:
we are interested in extending a given graph to include more vertices, while
forcing the extended graph to preserve the degree sequence of the original graph.


Let~\(G = (V(G),E(G))\) be a graph with vertex set~\(V(G) = \{1,2,\ldots,n\}\)
and edge set~\(E(G).\) Edges between vertices~\(i\) and~\(j\) are denoted as~\((i,j).\)
Two vertices~\(u\) and~\(v\) are said to be adjacent in~\(G\) if the edge~\((u,v) \in E(G)\)
and two edges are said to be adjacent if they share a common endvertex.
The degree~\(d_G(v)\) of a vertex~\(v\) denotes the number of vertices adjacent to~\(v\) in~\(G.\)
Let~\((d_G(1),\ldots,d_G(n))\) denote the degree sequence of the graph~\(G.\)
\begin{Definition}\label{def1} For integers~\(r, k \geq 1\) we say
that the graph~\(G\) is~\((r,k)-\)extendable or simply extendable
if there exists a graph~\(H\) with the following properties:\\
\((p1)\) The vertex set of~\(H\) is~\(\{1,2,\ldots,n+r\}.\)\\
\((p2)\) The degree sequence of~\(H\) is~\((d_G(1),\ldots,d_G(n),k,\ldots,k),\) where
the vertex~\(i\) has degree~\(d_G(i)\) for~\(1 \leq i \leq n\) and the vertex~\(n+j\)
has degree~\(k\) for~\(1 \leq j \leq r.\)
\end{Definition}
If there exists a graph~\(H\) satisfying~\((p1)-(p2),\)
then it is necessary that the product~\(r \cdot k\) is even,
since~\(\sum_{i=1}^{n}d_G(i)\) is even. We henceforth assume that~\(r \cdot k\) is even.
The graph~\(H\) in Definition~\ref{def1} is defined to
be an~\((r,k)-\)extension or simply, an extension of the graph~\(G.\)
For example, if~\(G\) is the graph consisting
of the single edge~\((1,2)\) with endvertices~\(\{1,2\},\)
then~\(G\) is \((1,2)-\)extendable and the graph~\(H_1\) with two edges~\((1,3)\) and~\((2,3)\)
is a~\((1,2)-\)extension of~\(G.\) The graph~\(G\)
is also \((2,1)-\)extendable and the graph~\(H_2\)
with two edges~\((1,3)\) and~\((2,4)\) is a~\((2,1)-\)extension of~\(G.\)


An edge of the extended graph~\(H\) with one endvertex in~\(\{n+1,\ldots,~n~+~r\}\)
and one endvertex in~\(\{1,2,\ldots,n\}\) is defined to be a \emph{cross} edge.
If~\(H\) contains no cross edges, we denote~\(H\) to be a \emph{trivial} extension of~\(G\) and~\(G\)
to be trivially~\((r,k)-\)extendable. For example, the single edge graph~\(G= \{(1,2)\}\)
is trivially~\((2,1)-\)extendable and the graph~\(H_3\) with two edges~\((1,2)\)
and~\((3,4)\) is a trivial~\((2,1)-\)extension of~\(G.\)

If~\(G\) is extendable but has no trivial extensions, then it is necessary to modify
the edge set and the vertex set of~\(G\) to obtain the desired extension~\(H.\)
To quantify the amount of modification needed, we introduce the concept
of edit number of graphs.

\subsection*{Edit number}
For integers~\(r,k \geq 1\) and a graph~\(G,\) the~\((r,k)-\)\emph{edit number} or simply edit number of~\(G\) is denoted by~\(N_{r,k}(G)\) and is defined as
follows. If~\(G\) is not~\((r,k)-\)extendable, then set~\(N_{r,k}(G) = \infty.\) Else set
\begin{equation}\label{nrk_def}
N_{r,k}(G) := \min_{H} \# (E(G) \setminus E_1(H) \cup E_1(H) \setminus E(G)),
\end{equation}
where the minimum is taken over all~\((r,k)-\)extensions~\(H\) of~\(G.\) The term\\\(E_1(H) \subset E(H)\) denotes the set of edges of~\(H\) with an endvertex in~\(\{1,2,\ldots,n\}.\) For example, if~\(G\) is the graph formed by the edge~\((1,2),\) then the only possible~\((1,2)-\)extension of~\(G\) is the graph~\(H_1\) with  edge set~\(\{(1,3), (2,3)\}\) and so~\(N_{1,2}(G) = 3.\) On the other hand, there are two possible~\((2,1)-\)extensions of~\(G:\) the graph~\(H_2 = \{(1,3), (2,4)\}\) and the trivial extension formed by the graph~\(H_3 = \{(1,2),(3,4)\}.\) For the trivial extension,~\(E(G) = E_1(H_3)\) and so~\(N_{2,1}(G) = 0.\)

In general for any~\(r,k \geq 1,\) the edit number~\(N_{r,k}(G) =0\) if and only if~\(G\) is trivially~\((r,k)-\)extendable. 
For all other extendable graphs~\(G,\) the edit number is nonzero and bounded above using~(\ref{nrk_def}) as
\begin{equation}\label{nrk_up}
N_{r,k}(G) \leq \#E(G) + \#E(H) = \frac{1}{2}\sum_{i=1}^{n} d_G(i) + \frac{1}{2}\sum_{i=1}^{n} d_G(i) + \frac{1}{2}r\cdot k.
\end{equation}


In what follows, we study graphs with low edit numbers based on the following interpretation of~\(N_{r,k}(G).\) Suppose~\(H\)
is a nontrivial~\((r,k)-\)extension of~\(G\) and we want to construct~\(H\)
from the graph~\(G \cup \{n+1,\ldots,n+r\}\) formed by the union of~\(G\)
and the isolated vertices~\(\{n+1,\ldots,n+r\}.\) This can be done in two steps as follows:
In the first step, add the edges of~\(E_1(H) \setminus E(G)\) and remove the edges of~\(E(G) \setminus E_1(H).\)
In the second step, add the edges of~\(E(H) \setminus E_1(H).\) If adding or removing an edge with
one endvertex in~\(\{1,2,\ldots,n\}\) constitutes an \emph{edit operation} on~\(G,\)
then from~(\ref{nrk_def}) we conclude that at least~\(N_{r,k}(G)\) edit operations are needed to obtain~\(H.\)
Moreover, there exists an extended graph~\(H_{opt}\) that can be obtained from~\(G\) after exactly~\(N_{r,k}(G)\) edit operations.

In Theorem~\ref{extn1} below, we classify all graphs with least possible edit number. For integers~\(p, q \geq 1,\) a subgraph~\(W\) of a graph~\(G\) is
said to be a~\((p,q)-\)subgraph of~\(G\) if it satisfies the following two properties:\\
\((i)\) The edge set~\(E(W)\) of~\(W\) contains~\(q\) edges and the vertex set~\(V(W)\) of~\(W\) is the set of all endvertices of~\(E(W).\)\\
\((ii)\) Every vertex of~\(W\) is adjacent to at most~\(p\) other vertices in~\(W.\)\\
Recall that if~\(G\) is~\((r,k)-\)extendable for some~\(r\) and~\(k,\) it is necessary that~\(r \cdot k\) is even
(see statement following Definition~\ref{def1}).
\begin{Theorem}\label{extn1}
Let~\(G\) be a graph with vertex set~\(\{1,2,\ldots,n\}\) and let~\(r, k \geq 1\) be integers
such that~\(r \cdot k\) is even.\\
\((1)\) If~\(1 \leq k \leq r-1,\) then~\(N_{r,k}(G) = 0.\)\\
\((2)\) If~\(k \geq r,\) then
\begin{equation}\label{nrg2}
N_{r,k}(G) \geq \frac{3}{2}r(k-r+1) 
\end{equation}
with equality if and only if~\(G\) contains an~\((r,t)-\)subgraph for~\(t = \frac{1}{2}r(k-r+1).\)
\end{Theorem}
In case~\((1),\) the graph~\(G\) has a trivial extension obtained from the union of~\(G\) and a regular graph~\(F\) on~\(r\) vertices, each having degree~\(k.\) We discuss construction of regular graphs in the next subsection.

A graph~\(G\) for which equality holds in~(\ref{nrg2}) is defined to be optimally\\\((r,k)-\)extendable or optimally extendable.
From Theorem~\ref{extn1}, we therefore get that \(G\) is optimally~\((r,k)-\)extendable if and only if~\(G\) contains an~\((r,t)-\)subgraph for~\(t = \frac{1}{2}r(k-r+1).\) The following Corollary gives sufficient conditions for a graph to be optimally extendable in
terms of matchings. A matching of size~\(l\) or a~\(l-\)matching of a graph~\(G\) is a set of~\(l\)
\emph{vertex disjoint} edges of~\(G.\) Two matchings~\(M_1\) and~\(M_2\) are said to be disjoint if~\(M_1\) and~\(M_2\) contain no edge in common.
\begin{Corollary}\label{extn} Let~\(G\) be a graph with vertex set~\(\{1,2,\ldots,n\}\) and let~\(r, k \geq 1\) be integers
such that~\(r \cdot k\) is even.\\
\((1)\) If~\(r=1\) and~\(k\) is even, then~\(G\) is optimally~\((1,k)-\)extendable if and only if~\(G\) contains a~\(\frac{k}{2}-\)matching.\\
\((2)\) If~\(k=r\) is even, then~\(G\) is optimally~\((k,k)-\)extendable if and only if~\(G\) contains~\(\frac{k}{2}\) edges.\\
\((3)\) If~\(k-r \geq 1\) is odd, then~\(G\)
is optimally~\((r,k)-\)extendable if~\(G\) contains~\(r\) disjoint~\(\left(\frac{k-r+1}{2}\right)-\)matchings.\\
\((4)\) If~\(k-r \geq 1\) is even, then~\(G\)
is optimally~\((r,k)-\)extendable if~\(G\) contains~\(r\) disjoint matchings, each of size~\(\frac{1}{2}(k-r)+1.\)\\
\((5)\) Suppose~\(d_G(i) \geq k\) for all~\(1 \leq i \leq n.\) If~\(k\) is even, then~\(G\) is optimally\\\((1,k)-\)extendable and if~\(k\) is odd,
then~\(G\) is optimally~\((2,k)-\)extendable.
\end{Corollary}





\subsection*{Regular graphs}
For integers~\(n,k \geq 1,\) a graph~\(G\) with vertex set~\(\{1,2,\ldots,n\}\) is said to be~\(k-\)\emph{regular} if the degree of each vertex in~\(G\)
is exactly~\(k.\) The graph~\(G\) is \emph{nearly}~\(k-\)regular if~\(n-1\) vertices of~\(G\) have degree~\(k\) and one vertex has degree~\(k-1.\) If~\(G\) is~\(k-\)regular, then the product~\(n\cdot k\) is the sum of degrees of the vertices in~\(G,\) which in turn is twice the number of edges in~\(G.\) Therefore it is necessary that either~\(n\) or~\(k\) is even.

In the proof of Theorem~\ref{thm1} below, we describe an algorithm to alternately construct connected regular graphs. A path~\(P = (v_1,\ldots,v_t)\) in the graph~\(G\) is a sequence of distinct vertices such that~\(v_i\) is adjacent to~\(v_{i+1}\) for~\(1 \leq i \leq t-1.\) The vertices~\(v_1\) and~\(v_t\) are connected by the path~\(P\) in~\(G.\) The graph~\(G\) is said to be \emph{connected} if any two vertices in~\(G\) are connected by a path in~\(G\) (Bollob\'as~(2002)).
\begin{Theorem}\label{thm1}
The following statements hold. \\
\((a)\) For all even integers~\(k \geq 2\) and all integers~\(n \geq k+1,\) there is a connected~\(k-\)regular graph with vertex set~\(\{1,2,\ldots,n\}.\)\\
\((b)\) For all odd integers~\(k \geq 3\) and all even integers~\(n \geq k+1,\) there is a connected~\(k-\)regular graph with vertex set~\(\{1,2,\ldots,n\}.\)\\
\((c)\) For all odd integers~\(k \geq 3\) and all odd integers~\(n \geq k+2,\) there is a connected nearly~\(k-\)regular graph with vertex set~\(\{1,2,\ldots,n\}.\)
\end{Theorem}
It is possible to use the Erd\H{o}s-Gallai Theorem (see Erd\H{o}s and Gallai (1960)) or the Havel-Hakimi criterion (Havel (1955), Hakimi (1962)) iteratively, to determine the existence of~\(k-\)regular graphs on~\(n\) vertices. In Theorem~\ref{thm1}, we provide a deterministic iterative algorithm to alternately construct \emph{connected}~\(k-\)regular graphs on~\(n\) vertices for all permissible values of~\(n \geq k+1,\) starting from the complete graph~\(K_{k+1}.\)

The paper is organized as follows. In Section~\ref{pf2}, we prove Theorem~\ref{thm1} and in Section~\ref{pf1}, we prove Theorem~\ref{extn1} and Corollary~\ref{extn}.

\setcounter{equation}{0}
\renewcommand\theequation{\thesection.\arabic{equation}}
\section{Proof of Theorem~\ref{thm1}} \label{pf2}
The following fact is used throughout.\\
\((a1)\) If the minimum degree of a vertex in a graph~\(G\) is~\(\delta \geq 2,\)
then there exists a path in~\(G\) containing~\(\delta\) edges.\\
\emph{Proof of~\((a1)\)}: Let~\(P = (v_1,\ldots,v_t)\) be the longest path in~\(G;\) i.e.,~\(P\) is a path containing the maximum number of edges.
Since~\(P\) is the longest path, all the neighbours of~\(v_1\) in~\(G\) belong to~\(P.\) Since~\(v_1\)
has at least~\(\delta\) neighbours in~\(G,\) we must have~\(t \geq \delta+1\) and so~\(P\) has at least~\(\delta\) edges.~\(\qed\)

We prove~\((a),(b)\) and~\((c)\) separately below.\\
\emph{Proof of~\((a)\)}: For~\(n = k+1,\) the complete graph~\(K_{n}\) is a connected~\(k-\)regular graph on~\(n\) vertices, with~\(k\) even. Suppose~\(G_{n}(k)\) is a connected~\(k-\)regular graph with vertex set~\(\{1,2,\ldots,n\},\) for some~\(n \geq k+1.\) Use property~\((a1)\) and let~\(P_n(k)\) be a path in~\(G_n(k)\) containing~\(k\) edges. Since~\(k\) is even, there are~\(\frac{k}{2}\) \emph{vertex disjoint} edges~\(e_i = (u_i,v_i), 1 \leq i \leq \frac{k}{2}\) in~\(P_n(k);\) i.e., the set~\(\bigcup_{1 \leq  i\leq \frac{k}{2}}\{u_i,v_i\}\) has~\(k\) distinct vertices. Remove the edges~\(e_i, 1 \leq i \leq \frac{k}{2}\) and add~\(k\) new edges~\[\bigcup_{i=1}^{\frac{k}{2}}\{(n+1,u_i)\}\bigcup \{(n+1,v_i)\}.\] The resulting graph is a connected~\(k-\)regular graph with vertex set\\\(\{1,2,\ldots,n+1\}\) and is defined to be~\(G_{n+1}(k).\)\\

\emph{Proof of~\((b)\)}: The proof is analogous as in the case of~\((a)\) above. For~\(n = k+1,\) the complete graph~\(K_{n}\) is a connected~\(k-\)regular graph on~\(n\) vertices with~\(k\) odd. Suppose~\(G_{n}(k)\) is a connected~\(k-\)regular graph with vertex set~\(\{1,2,\ldots,n\},\) for some even~\(n \geq k+1.\) We use~\(G_{n}(k)\) to construct a connected~\(k-\)regular graph~\(G_{n+2}(k)\) with vertex set~\(\{1,2,\ldots,n+2\}\) as follows.

By property~\((a1),\) the graph~\(G_n(k)\) contains a path~\(Q_n(k) = (q_1,\ldots,q_{k+1})\) consisting of~\(k\) edges~\(\{(q_j,q_{j+1})\}_{1 \leq  j \leq k}.\) Remove the~\(k-1\) edges~\(\{(q_j,q_{j+1})\}_{1 \leq j \leq k-1}\) and add the following edges:\\
\((i)\) For~\(1 \leq j \leq k-1, j\) odd, add the edges~\(\{(n+1,q_j),(n+1,q_{j+1})\}.\)\\
\((ii)\) For~\(1 \leq j \leq k-1, j\) even, add the edges~\(\{(n+2,q_j),(n+2,q_{j+1})\}.\)\\
\((iii)\) Add the edge~\((n+1,n+2).\)

Since~\(k\) is odd, the total number of edges added in step~\((i)\) is~\(k-1\) and so there are~\(k-1\) edges with~\(n+1\) as an endvertex after step~\((i).\)
Similarly, after step~\((ii)\) there are~\(k-1\) edges with~\(n+2\) as an endvertex. Finally, the resulting graph after step~\((iii)\) is a connected~\(k-\)regular graph with vertex set~\(\{1,2,\ldots,n+1,n+2\}\) and is defined to be~\(G_{n+2}(k).\)\\


\emph{Proof of~\((c):\)} Let~\(k\geq 2\) be odd and let~\(n \geq k+2\) be odd. We use the graph~\(G_{n-1}(k)\) with vertex set~\(\{1,2,\ldots,n-1\}\) obtained in~\((b)\) above to construct the graph~\(G_n(k)\) with vertex set~\(\{1,2,\ldots,n\}.\) From property~\((a1),\) the graph~\(G_{n-1}(k)\) contains a path~\(S_{n-1}(k)\) consisting of~\(k\) edges. Since~\(k\) is odd, there are~\(\frac{k-1}{2}\) vertex disjoint edges~\(f_i = (x_i,y_i), 1 \leq i \leq \frac{k-1}{2}\) in~\(S_{n-1}(k);\) i.e., the set~\(\bigcup_{1 \leq  i\leq \frac{k-1}{2}}\{x_i,y_i\}\) has~\(k-1\) distinct vertices. Remove the edges~\(f_i, 1 \leq i \leq \frac{k-1}{2}\) and add~\(k-1\) new edges~\[\bigcup_{i=1}^{\frac{k-1}{2}}\{(n,x_i)\}\bigcup \{(n,y_i)\}\] and define the resulting graph to be~\(G_{n}(k).\) By construction,~\(G_n(k)\) is connected, the vertex~\(n\) has degree~\(k-1\) and the rest of all the vertices have degree~\(k.\)~\(\qed\)



\setcounter{equation}{0}
\renewcommand\theequation{\thesection.\arabic{equation}}
\section{Proof of Theorem~\ref{extn1} and Corollary~\ref{extn}} \label{pf1}
We prove Theorem~\ref{extn1} first. To prove~\((1),\) Theorem~\ref{extn1} we use Theorem~\ref{thm1} and let~\(F\) be a~\(k-\)regular graph on~\(r\) vertices for~\(k \leq r-1.\) The union~\(F \bigcup G\) is a trivial~\((r,k)-\)extension of the graph~\(G\) and so~\(N_{r,k}(G) = 0.\) To prove~\((2)\) of Theorem~\ref{extn1}, we proceed in two steps. In the first step, we let~\(G\) be any graph such that~\(N_{r,k}(G) \leq 3t = \frac{3}{2}r(k-r+1)\) and using definition~(\ref{nrk_def}), we also let~\(Q\) be an~\((r,k)-\)extended graph such that
\begin{equation}\label{nrg_q}
\# (E(G) \setminus E_1(Q) \cup E_1(Q) \setminus E(G)) = N_{r,k}(G) \leq \frac{3}{2} r(k-r+1).
\end{equation}
We prove that such a graph~\(Q\) necessarily contains an~\((r,t)-\)subgraph and that~(\ref{nrg_q}) holds with equality. In the second step,
we let~\(G\) be any graph containing an~\((r,t)-\)subgraph~\(W\) for~\(t = \frac{1}{2} r(k-r+1).\) We construct an\\\((r,k)-\)extension of~\(G\) and show that~\(N_{r,k}(G)~=~3t,\) completing the proof of~\((2),\) Theorem~\ref{extn1}.

\underline{\emph{Step 1}}: We use an iterative colouring procedure on the edges of\\\(E_1(Q) \setminus E(G)\) and~\(E(G) \setminus E_1(Q),\)
briefly described as follows. In the first step of the iteration, we find a cross edge~\(e_1 \in E_1(Q) \setminus E(G)\) and colour it green. Recall (see paragraph following Definition~\ref{def1}) that cross edges have one endvertex in~\(\{n+1,\ldots,n+r\}\) and the other endvertex in~\(\{1,2,\ldots,n\}.\) In the next step, we find an edge~\(e_2 \in E(G) \setminus E_1(Q)\) adjacent to~\(e_1\) colour it red.
Such an edge must exist since the endvertex~\(v\) common to~\(e_1\) and~\(e_2\) has same degree in both the graphs~\(G\) and~\(Q.\)
In the third step, we again find an edge~\(e_3 \in E_1(Q) \setminus E(G)\) and colour it green. This procedure terminates
after a finite number of steps and in the final step, we pick a cross edge and colour it green. We then repeat the procedure picking a new (uncoloured)
cross edge and continue this way until all cross edges are exhausted. We then use these coloured edges and relation~(\ref{nrg_q})
to obtain the desired~\((r,t)-\)subgraph.

For a vertex~\(v \in \{1,2,\ldots,n\},\) let~\(f^{(0)}_Q(v)\) and~\(f^{(0)}_G(v)\) respectively denote the number of edges in~\(E_1(Q) \setminus E(G)\) and~\(E(G)\setminus E_1(Q)\) containing~\(v\) as an endvertex. Since the degrees of~\(v\) in the graphs~\(G\) and~\(Q\) are the same, i.e.,~\(d_G(v) = d_Q(v),\) we must have that~\(f^{(0)}_Q(v) = f^{(0)}_G(v)\) for all~\(v \in \{1,2,\ldots,n\}.\)


Set~\(i=1\) and~\(x_0 = n+1.\)
Pick an edge~\((x_{i-1},x_i) \in E_1(Q)\setminus E(G)\) and colour it green. The term~\(f_Q^{(i-1)}(x_i) \geq 1\) necessarily and we set
\begin{equation}\label{fiq}
f^{(i)}_Q(x_i) := f^{(i-1)}_Q(x_i)-1, \;\;\;f^{(i)}_G(x_i) := f^{(i-1)}_G(x_i)
\end{equation}
and for all~\(v \in \{1,2,\ldots,n\}, v \neq x_i,\) we set
\begin{equation}\label{fiq2}
f^{(i)}_Q(v) := f^{(i-1)}_Q(v) \text{ and }f^{(i)}_G(v) := f^{(i-1)}_G(v).
\end{equation}
The terms~\(f^{(i)}_Q(v)\) and~\(f^{(i)}_G(v)\) denote, respectively, the number of remaining uncoloured edges of~\(E_1(Q) \setminus E(G)\) and~\(E(G) \setminus E_1(Q)\) containing~\(v\) as an endvertex, after iteration step~\(i.\) Since~\(f^{(i-1)}_Q(x_i) = f^{(i-1)}_G(x_i),\) we get from~(\ref{fiq}) that
\begin{equation}\label{fiq3}
f^{(i)}_G(x_i) = f^{(i)}_Q(x_i)+1.
\end{equation}
This completes the first step of the iteration.

In the next step, we set~\(i=2\) and obtain from~(\ref{fiq3}) that~\(f^{(i-1)}_G(x_{i-1}) \geq 1\) and so there exists an uncoloured edge~\((x_{i-1},x_i) \in E(G) \setminus E_1(Q)\) containing~\(x_{i-1}\) as an endvertex. Moreover, the vertex~\(x_i \in \{1,2,\ldots,n\}\) and\\~\(f^{(i-1)}_G(x_i) \geq 1.\) We colour the edge~\((x_{i-1},x_i)\) red and set
\begin{equation}\label{fiqq}
f^{(i)}_G(x_i) := f^{(i-1)}_G(x_i)-1, \;\;\;f^{(i)}_Q(x_i) := f^{(i-1)}_Q(x_i)
\end{equation}
and for all~\(v \in \{1,2,\ldots,n\}, v \neq x_i,\) we set
\begin{equation}\label{fiqq2}
f^{(i)}_Q(v) := f^{(i-1)}_Q(v) \text{ and }f^{(i)}_G(v) := f^{(i-1)}_G(v).
\end{equation}
We recall from the previous update~(\ref{fiq2}) that~\(f^{(i-1)}_Q(x_i) = f^{(i-1)}_G(x_i)\) and so from~(\ref{fiqq}) we get that
\begin{equation}\label{fiqq3}
f^{(i)}_Q(x_i) = f^{(i)}_G(x_i)+1.
\end{equation}
This completes the second step of the iteration.

In the third step we set~\(i=3\) and arguing as before, obtain an edge\\\((x_{i-1},x_i) \in E_1(Q) \setminus E(G)\) and colour it green. If~\(x_i \notin \{n+1,\ldots, n+r\}\) we then find an edge of~\(E(G) \setminus E_1(Q)\) containing~\(x_i\) as an endvertex and colour it red.
Since the number of vertices in~\(G\) is finite, the above procedure of alternately colouring edges green and red terminates at some iteration step~\(j,\) after we reach a vertex~\(x_{j} \in \{n+1,\ldots,n+r\}.\) Moreover, the edge~\((x_{j-1},x_{j}) \in E_1(Q) \setminus E(G)\) is coloured green
and
\begin{equation}\label{f_eq_fin}
f^{(j)}_G(v) = f^{(j)}_Q(v)
\end{equation}
for all~\(v \in \{1,2,\ldots,n\}.\)


Let
\begin{equation}\label{cj_def}
C_1(1) = (n+1,x_1,x_2,\ldots,x_{j-1},x_{j})
\end{equation}
denote the resulting sequence of edges obtained after the above iteration terminates. We denote~\(C_1(1)\) to be a \emph{walk}. Let~\(P_1(1)\) be the subwalk of~\(C_1(1)\) with endvertices~\(x_1\) and~\(x_{j-1}\) and not containing the vertices~\(n+1\) or~\(x_{j}.\)

Retain the colouring of the edges obtained above and relabel~\(f^{(j)}_Q(v)\) and~\(f^{(j)}_G(v)\) as~\(f^{(0)}_Q(v)\) and~\(f^{(0)}_G(v),\) respectively.
From~(\ref{f_eq_fin}), we have that\\\(f^{(0)}_Q(v) = f^{(0)}_G(v)\) for all~\(v \in \{1,2,\ldots,n\}.\) We then repeat the above procedure starting with an uncoloured edge containing~\(n+1\) as an endvertex and obtain another walk~\(C_1(2)\) and a corresponding subwalk~\(P_1(2)\) with endvertices in~\(T.\) Continue this procedure until all edges with~\(n+1\) as an endvertex are coloured and suppose we obtain~\(\alpha_1\) walks~\(\{C_1(j)\}_{1 \leq j \leq \alpha_1}\) and corresponding subwalks~\(\{P_1(j)\}_{1 \leq j \leq \alpha_1}.\)

Retain the colouring and perform the same procedure as above with the remaining uncoloured edges containing~\(n+2\) as an endvertex to obtain~\(\alpha_2\) walks~\(\{C_2(j)\}_{1 \leq j \leq \alpha_2}\) and corresponding subwalks~\(\{P_2(j)\}_{1 \leq j \leq \alpha_2}\) having endvertices in~\(T.\) Continue
this colouring procedure until all edges having an endvertex in~\(\{n+1,\ldots,n+r\}\) are coloured and let~\(\{C_i(j)\}_{1 \leq i \leq r, 1 \leq j \leq \alpha_i}\) and~\(\{P_i(j)\}_{1 \leq i \leq r, 1 \leq j \leq \alpha_i}\) be the set of all walks obtained. The following properties hold:\\
\((w3)\) All the green edges of~\(C_i(j)\) belong to~\(E_1(Q) \setminus E(G)\) and all the red edges of~\(C_i(j)\) belong to~\(E(G) \setminus E_1(Q).\)
Each subwalk~\(P_i(j)\) contains at least one (red) edge.\\
\((w4)\) The total number of walks is
\begin{equation}\label{alp_sum}
\sum_{i=1}^{r} \alpha_i \geq \frac{1}{2}r(k-r+1) = t.
\end{equation}
Property~\((w3)\) is true by construction. To prove~\((w4),\) define any edge of the graph~\(Q\) containing one endvertex in~\(\{n+1,\ldots,n+r\}\) and other endvertex in~\(\{1,2,\ldots,n\}\) to be a cross edge. Since~\(v \in \{n+1,\ldots,n+r\}\) is adjacent to at most~\(r-1\) other vertices of~\(\{n+1,\ldots,n+r\},\) the vertex~\(v\) is adjacent to at least~\(k-r+1\) vertices of~\(G.\) Therefore, the total number of cross edges in~\(Q\) is atleast~\(r(k-r+1).\) From the construction above, every walk in~\(\{C_i(j)\}_{1 \leq i \leq r, 1 \leq j \leq \alpha_i}\) contains two distinct cross edges and the walks~\(\{C_i(j)\}\) are edge disjoint. Thus~(\ref{alp_sum}) holds.

Denoting~\(\#C_i(j)\) and~\(\#P_i(j)\) to be the number of edges in the walks~\(C_i(j)\) and~\(P_i(j),\) respectively, we get from property~\((w3)\) that
\begin{eqnarray}
 \sum_{i=1}^{r} \sum_{j=1}^{\alpha_i} \left(\#P_i(j)+2\right)  &=& \sum_{i=1}^{r} \sum_{j=1}^{\alpha_i} \#C_i(j)\nonumber\\
&\leq& \#(E_1(Q) \setminus E(G)) + \#(E(G) \setminus E_1(Q)) \nonumber\\
&\leq& \frac{3}{2}r(k-r+1). \label{eq_up}
\end{eqnarray}
The final estimate in~(\ref{eq_up}) follows from~(\ref{nrg_q}).



Using~(\ref{alp_sum}) in~(\ref{eq_up}), we get
\begin{equation} \label{alp_sum4}
0 \leq \sum_{i=1}^{r} \sum_{j=1}^{\alpha_i} \left(\#P_i(j) -1\right) \leq 0
\end{equation}
and so from~\((w3),\) each subwalk~\(P_{i}(j)\) consists of a single (red) edge. This in turn implies that the walk~\(C_i(j)\) contains exactly two green cross edges and the red edge of~\(P_{i}(j).\) Also, the graph~\(Q\) is constructed in the following two steps:\\
~\((s1)\) Remove all the (red) edges of~\(\{P_i(j)\}\) and add all the (green) edges of~\(\{C_i(j)\setminus P_i(j)\}.\)\\
~\((s2)\) Add all the~\({r \choose 2}\) edges with both endvertices in~\(\{n+1,\ldots,n+r\}.\)\\
From~\((s1)-(s2)\) we then get
\begin{equation}\label{ep_we}
\bigcup_{i,j} P_i(j) = E(G) \setminus E_1(Q) \text{ and } \bigcup_{i,j} (C_i(j) \setminus P_i(j)) = E_1(Q) \setminus E(G).
\end{equation}

To estimate the total number of edges in the walks~\(\{C_i(j)\}\) and~\(\{P_i(j)\},\) substitute~\(\#P_i(j) =1\) in the first term of~(\ref{eq_up}) to get~\[\sum_{i=1}^{r}\alpha_i \leq \frac{1}{2}r(k-r+1).\] From~(\ref{alp_sum}), we then conclude that
\begin{equation}\label{epij}
\sum_{i=1}^{r} \alpha_i = \frac{1}{2}r(k-r+1) = \sum_{i=1}^{r} \sum_{j=1}^{\alpha_i} \#P_i(j)
\end{equation}
and so the number of edges in~\(\{P_i(j)\}\) is~\(t = \frac{1}{2}r(k-r+1).\) Recall that by construction, each walk~\(C_i(j)\)
contains two green edges and the walks~\(\{C_i(j)\}\) are all edge disjoint. So~\[\sum_{i=1}^{r} \sum_{j=1}^{\alpha_i} \#(C_i(j) \setminus P_i(j))= 2\sum_{i=1}^{r}\alpha_i = 2t,\] using~(\ref{epij}).
From~(\ref{ep_we}) and~(\ref{nrg_q}) we therefore get that
\begin{equation}\label{nr_eqw}
N_{r,k}(G) = \#(E(G) \setminus E_1(Q)) + \#(E_1(Q)\setminus E(G)) = 3t.
\end{equation}


Finally, to obtain an~\((r,t)-\)subgraph of~\(G,\) let~\(W_0 \subseteq G\) be the graph formed by the edges~\(\bigcup_{i,j} P_{i}(j) = E(G) \setminus E_1(Q)\)
and suppose vertex~\(v\) has degree~\(d\) in~\(W_0.\) Let~\(W_0(v) \subset W_0 \) be the~\(d\) edges containing~\(v\) an endvertex and let~\(Y_0(v) \subseteq Q\) be the set of edges in~\(E_1(Q) \setminus E(G)\) containing~\(v\) as an endvertex.  Since the degree of~\(v\) is the same in~\(G\) and~\(Q\) we must have~\(\#Y_0(v) = \#W_0(v) = d.\) Also, from discussion prior to~(\ref{ep_we}), every edge in~\(Y_0(v)\) has an endvertex in~\(\{n+1,\ldots,n+r\}\) and so~\(v\) is adjacent to exactly~\(d\) vertices of~\(\{n+1,\ldots,n+r\}\) in the graph~\(Q.\) Thus~\(d \leq r\) and the graph~\(W_0\) is the desired~\((r,t)-\)subgraph of~\(G.\) This completes the proof of the statement following~(\ref{nrg_q}).

\underline{\emph{Step 2}}: Suppose the graph~\(G\) contains an~\((r,t)-\)subgraph~\(W\) consisting of~\(t = \frac{1}{2} r(k-r+1)\) edges. We construct an~\((r,k)-\)extension of~\(G\) and obtain that~\(N_{r,k}(G)~=~3t.\) Let~\(\{e_1,\ldots,e_t\}\) and~\(\{v_1,\ldots,v_x\}\) respectively denote the edge set and vertex set of the graph~\(W.\) If~\(v_i, 1 \leq i \leq x\) is adjacent to~\(l_i \leq r\) vertices of~\(W,\) define~\(v_i(j) = v_i\) for~\(1 \leq j \leq l_i.\) The vector~\[V = (v_1(1),\ldots,v_1(l_1),v_2(1),\ldots,v_2(l_2),\ldots,v_x(1),\ldots,v_x(l_x))\]
has~\(2t\) entries, which we relabel as~\(v_1(1) = z_1, v_1(2) = z_2\) and so on, so that~\(V = (z_1,\ldots,z_{2t}).\)

Remove the~\(t\) edges~\(\{e_1,\ldots,e_t\}\) of~\(W\) from the graph~\(G\) and add the~\({r \choose 2}\) edges
between vertices in~\(\{n+1,\ldots,n+r\}\) and call the resulting graph as~\(H_{temp}.\) For the case~\(k-r+1\) is even, perform the following additional connection procedure on the graph~\(H_{temp}.\) In the first step, connect the vertices~\(z_i\) and~\(z_{r+i}\) to the vertex~\(n + i,\) for~\(1 \leq i \leq r.\) In the second step, connect the vertices~\(z_{2r+i}\) and~\(z_{3r+i}\) to~\(n+i\) for~\(1 \leq i \leq r\) and so on. Thus for example, the vertex~\(n+1\) would be connected to the vertices~\(z_1,z_{r+1},z_{2r+1}\) and~\(z_{3r+1}\) after the first two steps. Since each vertex of the graph~\(W\) appears at most~\(r\) times in the vector~\(V,\) the vertices~\(z_1,z_{r+1},z_{2r+1}\) and~\(z_{3r+1}\) are all distinct.

The above procedure terminates after exactly~\(\frac{2t}{2r} = \frac{1}{2}(k-r+1)\) steps and in the resulting final
graph~\(H_{fin},\) each vertex in~\(\{n+1,\ldots,n+r\}\) has degree exactly~\(k.\) Also the degree of each vertex with~\(\{1,2,\ldots,n\}\)
is the same as the original graph~\(G.\) This is true because if a vertex~\(v\) has degree~\(l_v\) in~\(W,\) then~\(l_v\)
edges of~\(G\) adjacent to~\(v\) have been removed in obtaining~\(H_{temp}\) and by the connection procedure above,~\(l_v\) new cross edges with one endvertex as~\(v\) have been added. Thus~\(H_{fin}\) is the desired~\((r,k)-\)extension of~\(G.\)

If~\(k-r+1\) is odd, then~\(r\) is even since the product~\(r\cdot k\) is even and we perform the connection procedure in a slightly different way after obtaining the graph~\(H_{temp}.\) In the first step of the connection procedure, we connect vertices~\(z_i\) and~\(z_{r + i}\) to vertex~\(n+i\) for~\(1 \leq i \leq \frac{r}{2}.\) In the second step of the procedure, we connect vertices~\(z_{\frac{r}{2} + i}\) and vertices~\(z_{\frac{3r}{2} + i}\) to vertex~\(n+i\)
for~\(\frac{r}{2} \leq i \leq r.\) In the third step of the procedure, we again connect vertices~\(z_{r+i}\) and~\(z_{2r+i}\)
to vertex~\(n+i\) for~\(1 \leq i \leq \frac{r}{2}.\) This procedure is continued for~\(k-r\) steps after which~\(r\)
vertices~\(z_{2t-r+1},\ldots,z_{2t}\) are left out. Connect vertex~\(z_{2t-r+i}\) to vertex~\(n+i\) for~\(1 \leq i \leq r.\)
As before, we obtain the desired extended graph~\(H_{fin}.\)

From the construction of the extended graph~\(H_{fin}\) described above, we estimate the edit number of the original graph~\(G\) as follows. First, to obtain the graph~\(H_{temp}\) from~\(G\) we removed~\(t\) edges with an endvertex in~\(G.\) Since no edge with an endvertex in~\(G\) was removed in obtaining~\(H_{fin}\) from~\(H_{temp},\) we have
\begin{equation}\label{eg1}
\#E(G)\setminus E_1(H_{fin}) = t.
\end{equation}
In the connection procedure used to obtain~\(H_{fin}\) from~\(H_{temp},\) we added~\(2t\) edges to~\(H_{temp},\) each containing an endvertex in~\(G.\) No edge with an endvertex in~\(G\) was added in obtaining~\(H_{temp}\) from~\(G\) and so
\begin{equation}\label{eg2}
\#E_1(H_{fin}) \setminus E(G) = 2t.
\end{equation}
From~(\ref{nrk_def}),~(\ref{eg1}) and~(\ref{eg2}), we therefore get that
\begin{equation}\label{nrk_up1}
N_{r,k}(G) \leq 3t = \frac{3}{2}r(k-r+1).
\end{equation}
But from the proof of Step~1 (see statement following~(\ref{nrg_q})), we then get that~\(N_{r,k}(G) = 3t.\)~\(\qed\)

\emph{Proof of Corollary~\ref{extn}}:
The proof of~\((1)\) follows from the fact that any\\\((1,t)-\)subgraph of~\(G\) with~\(t = \frac{k}{2}\) is a~\(\frac{k}{2}-\)matching and vice versa. Property~\((2)\) is true since if~\(W\) is any collection of~\(\frac{k}{2}\) edges of~\(G,\) then each vertex in~\(W\) is adjacent to at most~\(k\) edges of~\(W\) and so~\(W\) is a~\((k,t)-\)subgraph of~\(G\) with~\(t = \frac{k}{2}.\)

We prove property~\((3)\) and the proof of~\((4)\) is analogous. To prove property~\((3),\) let~\(M_1,\ldots,M_r\) be~\(r\) disjoint~\(\left(\frac{k-r+1}{2}\right)-\)matchings of~\(G.\) It suffices to see that the collection of edges in~\(\{M_1,\ldots,M_r\}\) form an~\((r,t)-\)subgraph with~\(t =\frac{1}{2}r(k-r+1).\) The total number of edges in~\(\{M_1,\ldots,M_r\}\)
is~\(t.\) Also if~\(v\) is an endvertex of some edge in~\(\{M_1,\ldots,M_r\},\) then~\(v\) is adjacent to at most~\(r\) edges in~\(\{M_1,\ldots,M_r\},\) because otherwise two edges of some matching would both contain~\(v\) as an endvertex, a contradiction.

To prove property~\((5),\) we use property~\((a1)\) in the proof of Theorem~\ref{thm1} and obtain a path~\(P\) of length~\(k\) in~\(G.\)
If~\(k\) is even, then the alternate edges of~\(P\) form a~\(\frac{k}{2}-\)matching in~\(G\) and so~\(G\) is optimally~\((1,k)-\)extendable
by property~\((1).\) If~\(k\) is odd, the alternate edges of the first~\(k-1\) edges of~\(P\) and the remaining~\(\frac{k-1}{2}\) edges
form two disjoint~\(\left(\frac{k-1}{2}\right)-\)matchings in~\(G\) and so the result follows from property~\((3).\)~\(\qed\)

\subsection*{Acknowledgement}
I thank Professors Rahul Roy, Thomas Mountford and Federico Camia for crucial comments and for my fellowships.

\bibliographystyle{plain}

\end{document}